\documentclass[12pt,reqno]{amsart}
\usepackage{enumerate, latexsym, amsmath, amsfonts, amssymb, amsthm, color}
\def\Z{\Bbb Z}
\def\N{\Bbb N}

\def\l{\left}
\def\r{\right}

\def\t{\text}
\def\f{\frac}
\def\mo{{\rm{mod}\ }}

\def\ls{\leqslant}
\def\gs{\geqslant}

\def\bi{\binom}
\def\al{\alpha}

\def\eq{\equiv}

\def\em{\emptyset}
\def\Pr{{\rm Pr}}
\def\PR{{\rm PR}}

\def\Ack{\medskip\noindent {\bf Acknowledgments}}
\theoremstyle{plain}

\newtheorem{conjecture}{Conjecture}
\theoremstyle{definition}

\theoremstyle{remark}
\newtheorem{remark}{Remark}

 \vspace{4mm}

\begin{document}
 \baselineskip=16pt
\hbox{Number Theory: Arithmetic in Shangri-La (eds., S. Kanemitsu,
H. Li and J. Liu),} \hbox{Proc. 6th China-Japan Seminar (Shanghai,
August 15-17, 2011),} \hbox{World Sci., Singapore, 2013, pp.
244-258.}
\medskip

\title
[Conjectures involving arithmetical sequences]
{Conjectures involving \\arithmetical sequences}

\author
[Zhi-Wei Sun] {Zhi-Wei Sun}

\thanks{Supported by the National Natural Science Foundation (grant 11171140)
 of China}

\address {Department of Mathematics, Nanjing
University, Nanjing 210093, People's Republic of China}
\email{zwsun@nju.edu.cn}


 \begin{abstract} We pose thirty conjectures on arithmetical sequences, most
of which are about monotonicity of sequences of the form $(\root
n\of{a_n})_{n\gs 1}$ or the form $(\root{n+1}\of{a_{n+1}}/\root
n\of{a_n})_{n\gs1}$, where $(a_n)_{n\gs 1}$ is a number-theoretic or
combinatorial sequence of positive integers. This material might
stimulate further research.
\end{abstract}

\maketitle

\section{Introduction}
\setcounter{lemma}{0}
\setcounter{theorem}{0}
\setcounter{corollary}{0}
\setcounter{remark}{0}
\setcounter{equation}{0}
\setcounter{conjecture}{0}

 A sequence $(a_n)_{n\gs0}$ of natural numbers is said to be {\it
log-concave} (resp. {\it log-convex}) if $a_{n+1}^2\gs a_na_{n+2}$
(resp. $a_{n+1}^2\ls a_na_{n+2}$) for all $n=0,1,2,\ldots$. The
log-concavity or log-convexity of combinatorial sequences has been
studied extensively by many authors (see, e.g., \cite{C,D05,D09,J,LW}).

For $n\in\Z^+=\{1,2,3,\ldots\}$ let $p_n$ denote the $n$-th prime. In
1982, Faride Firoozbakht conjectured that
$$\root n\of{p_n}>\root{n+1}\of {p_{n+1}}\quad\t{for all}\ n\in\Z^+,$$
i.e., the sequence $(\root{n}\of {p_n})_{n\gs1}$ is strictly
decreasing (cf. \cite[p.\,185]{R}). This was verified for $n$ up to $3.495\times10^{16}$ by Mark Wolf \cite{W}.

Mandl's inequality (cf. \cite{Du,RS,H}) asserts that $S_n<np_n/2$ for all $n\gs 9$, where $S_n$ is the sum of the first $n$ primes.
Recently the author \cite{S12c} proved that
the sequence $(\root{n}\of {S_n})_{n\gs2}$ is strictly decreasing
and moreover the sequence $(\root{n+1}\of{S_{n+1}}/\root{n}\of
{S_n})_{n\gs5}$ is strictly increasing. Motivated by this, here we pose many
conjectures on sequences $(\root{n}\of{a_n})_{n\gs1}$ and
$(\root{n+1}\of{a_{n+1}}/\root{n}\of{a_n})_{n\gs1}$ for many
number-theoretic or combinatorial sequences $(a_n)_{n\gs1}$ of positive integers.
Clearly, if $(\root{n+1}\of{a_{n+1}}/\root{n}\of{a_n})_{n\gs N}$ is strictly increasing (decreasing) with limit 1, then
the sequence $(\root{n}\of{a_n})_{n\gs N}$ is strictly decreasing (resp., increasing).

Sections 2 and 3 are devoted to our conjectures involving number-theoretic sequences and combinatorial sequences respectively.

\section{Conjectures on number-theoretic sequences}
\setcounter{lemma}{0}
\setcounter{theorem}{0}
\setcounter{corollary}{0}
\setcounter{remark}{0}
\setcounter{equation}{0}
\setcounter{conjecture}{0}

\subsection{Conjectures on sequences involving primes}

\begin{conjecture}\label{Conj2.1} {\rm (2012-09-12)} For any $\al>0$ we have
$$\f1n\sum_{k=1}^n p_k^{\alpha}<\f{p_n^{\al}}{\al+1}\quad \t{for all}\ \ n\gs 2\lceil \al\rceil^2+\lceil \al\rceil+6.$$
\end{conjecture}

\begin{remark}\label{Rem2.1} We have verified the conjecture for $\al=2,3,\ldots,700$ and $n\ls 10^6$.
Our numerical computation suggests that for $\al=2,3,\ldots,10$ we may replace $\lceil\al\rceil^2+\lceil \al\rceil+6$ in the inequality by
$9,\, 15,\,31,\,47,\,62,\,92,\,92,\,122,\,122$
respectively. Note that Mandl's inequality (corresponding to the case $\al=1$)
can be restated as $\sum_{k=1}^n p_k<\f{n-1}2p_{n+1}$ for $n\gs8$, which provides a lower bound for $p_{n+1}$
in terms of $p_1,\ldots,p_n$.
\end{remark}

Our next conjecture is a refinement of Firoozbakht's conjecture.

\begin{conjecture}\label{Conj2.2} {\rm (2002-09-11)} For any integer $n>4$, we have the inequality
$$\f{\root{n+1}\of{p_{n+1}}}{\root n\of{p_n}}<1-\f{\log\log n}{2n^2}.$$
\end{conjecture}

\begin{remark}\label{Rem2.2} The author has verified the conjecture for all $n\ls3500000$ and all those $n$
with $p_n<4\times 10^{18}$ and $p_{n+1}-p_n\not=p_{k+1}-p_k$ for all $1\ls k<n$. Note that
if $n=49749629143526$ then $p_n=1693182318746371$, $p_{n+1}-p_n=1132$ and
$(1-\root{n+1}\of{p_{n+1}}/\root n\of{p_n})n^2/\log\log n\approx 0.5229.$
\end{remark}

A well-known theorem of Dirichlet (cf. \cite[pp.\,249-268]{IR}) states that for any relatively
prime positive integers $a$ and $q$ the arithmetic progression
$a,a+q,a+2q,\ldots$ contains infinitely many primes; we use $p_n(a,q)$ to denote the $n$-th prime in this progression.

The following conjecture extends the Firoozbakht conjecture to
primes in arithmetic progressions.

\begin{conjecture}\label{Conj2.3} {\rm (2012-08-11)} Let $q\gs a\gs1$  be
positive integers with $a$ odd, $q$ even and $\gcd(a,q)=1$.
Then there is a positive integer $n_0(a,q)$ such that the sequence $(\root{n}\of
{p_n(a,q)})_{n\gs n_0(a,q)}$ is strictly decreasing. Moreover, we may take $n_0(a,q)=2$ for $q\ls 45$.
\end{conjecture}

\begin{remark}\label{Rem2.3} Note that
$\root4\of{p_4(13,46)}<\root 5\of{p_5(13,46)}$. Also,  $\root3\of{p_3(3,328)}<\root4\of{p_4(3,328)}$
and $\root6\of{p_6(23,346)}<\root7\of{p_7(23,346)}$.
\end{remark}

A famous conjecture of E. Artin asserts that if $a\in\Z$ is neither $-1$ nor a square then there are infinitely many primes $p$ having $a$ as a primitive root
modulo $p$. This is still open, the reader may consult the survey \cite{M} for known progress on this conjecture.

\begin{conjecture}\label{Conj2.4} {\rm (2012-08-17)}  Let $a\in\Z$ be not a perfect power (i.e., there are no integers $m>1$ and $x$ with $x^m=a$).

{\rm (i)} Assume that $a>0$. Then there are infinitely many primes $p$ having $a$ as the smallest positive primitive root modulo $p$.
Moreover, if $p_1(a),\ldots,p_n(a)$ are the first $n$ such primes, then the next such prime $p_{n+1}(a)$ is smaller than $p_n(a)^{1+1/n}$, i.e.,
$\root n\of{p_n(a)}>\root{n+1}\of{p_{n+1}(a)}$.

{\rm (ii)} Suppose that $a<0$. Then there are infinitely many primes $p$ having $a$ as the largest negative primitive root modulo $p$.
Moreover, if $p_1(a),\ldots,p_n(a)$ are the first $n$ such primes, then the next such prime $p_{n+1}(a)$ is smaller than $p_n(a)^{1+1/n}$ (i.e.,
$\root n\of{p_n(a)}>\root{n+1}\of{p_{n+1}(a)}$) with the only exception $a=-2$ and $n=13$.

{\rm (iii)} The sequence $(\root{n+1}\of {P_{n+1}(a)}/\root n\of{P_n(a)})_{n\gs3}$ is strictly increasing with limit $1$,
where $P_n(a)=\sum_{k=1}^np_k(a)$.
\end{conjecture}

\begin{remark}\label{Rem2.4}  Let us look at two examples. The first 5 primes having 24 as the smallest positive
primitive root are $p_1(24)=533821$, $p_2(24)=567631$, $p_3(24)=672181$, $p_4(24)=843781$ and $p_5(24)=1035301$,
and we can easily verify that
$$p_1(24)>\sqrt{p_2(24)}>\root{3}\of {p_3(24)}>\root 4\of{p_4(24)}>\root 5\of{p_5(24)}.$$
The first prime having $-12$ as the largest negative primitive root is $p_1(-12)$ $=7841$, and the second prime having $-12$ as the largest negative
primitive root is $p_2(-12)=16061$; it is clear that $p_1(-12)>\sqrt{p_2(-12)}$.
\end{remark}

Recall that the Proth numbers have the form $k\times2^n+1$ with $k$ odd and $0<k<2^n$. In 1878 F. Proth proved that a Proth number $p$ is a prime
if (and only if) $a^{(p-1)/2}\eq-1\ (\mo\ p)$ for some integer $a$ (cf. Ex. 4.10 of \cite[p.\,220]{CP}).
A Proth prime is a Proth number which is also a prime number; the Fermat primes are a special kind of Proth primes.

\begin{conjecture}\label{Conj2.5} {\rm (2012-09-07)}
{\rm (i)} The number of Proth primes not exceeding a large integer $x$ is asymptotically equivalent to $c\sqrt x/\log x$
for a suitable constant $c\in(3,4)$.

{\rm (ii)} If $\Pr(1),\ldots,\Pr(n)$ are the first $n$ Proth primes, then the next Proth prime $\Pr(n+1)$ is smaller than $\Pr(n)^{1+1/n}$
 $($i.e., $\root n\of{\Pr(n)}>\root{n+1}\of{\Pr(n+1)})$ unless $n=2,4,5$. If we set $\PR(n)=\sum_{k=1}^n\Pr(k)$, then
  $\PR(n)<n\Pr(n)/3$ for all $n>50$, and
the sequence $(\root{n+1}\of{\PR(n+1)}/\root n\of{\PR(n)})_{n\gs 34}$
 is strictly increasing with limit $1$.
\end{conjecture}

\begin{remark}\label{Rem2.5} We have verified that $\root n\of{\Pr(n)}>\root{n+1}\of{\Pr(n+1)}$ for all $n=6,\ldots,4000$,
 $\PR(n)<n\Pr(n)/3$ for all $n=51,\ldots,3500$, and
$$\root{n+1}\of{\PR(n+1)}/\root n\of{\PR(n)}<\root{n+2}\of{\PR(n+2)}/\root {n+1}\of{\PR(n+1)}$$
 for all $n=34,\ldots,3200$.
\end{remark}

In the remaining part of this section,
 we usually list certain primes of special types in ascending order as $q_1,q_2,q_3,\ldots$, and write $Q(n)$ for $\sum_{k=1}^nq_k$.
 Note that the inequality $\root{n}\of{Q(n)}/\root {n-1}\of{Q(n-1)}<\root{n+1}\of{Q(n+1)}/\root n\of{Q(n)}$ yields a lower bound for $q_{n+1}$.

\begin{conjecture}\label{Conj2.6} {\rm (i) (2012-08-18)} Let $q_1,q_2,q_3,\ldots$ be the list (in ascending order)
of those primes of the form $x^2+1$ with $x\in\Z$.
Then we have $q_{n+1}<q_n^{1+1/n}$ unless $n=1,2,4,351$.
Also, the sequence $(\root{n+1}\of{Q(n+1)}/\root n\of{Q(n)})_{n\gs13}$ is strictly increasing with
limit $1$.

{\rm (ii) (2012-09-07)} Let $q_1,q_2,q_3,\ldots$ be the list (in ascending order) of those primes of the form $x^2+x+1$ with $x\in\Z$.
Then we have $q_{n+1}<q_n^{1+1/n}$ unless $n=3,6$. Also, the sequence $(\root{n+1}\of{Q(n+1)}/\root n\of{Q(n)})_{n\gs20}$ is strictly increasing with
limit $1$.
\end{conjecture}

\begin{remark}\label{Rem2.6} If we use the notation in part (i) of Conj. 2.6, then $q_{351}=3536^2+1=12503297$, $q_{352}=3624^2+1=13133377$, and
$\root{351}\of {q_{351}}<\root{352}\of{q_{352}}$.
\end{remark}

Schinzel's Hypothesis H (cf. \cite[pp.\,17-18]{CP}) states that if $f_1(x),\ldots,f_k(x)$ are irreducible polynomials
with integer coefficients and positive leading coefficients
such that there is no prime dividing the product $f_1(q)\cdots f_k(q)$ for all $q\in\Z$, then there are infinitely many $n\in\Z^+$
 such that $f_1(n),\ldots,f_k(n)$ are all primes.

 Here is a general conjecture related to Hypothesis H.

\begin{conjecture}\label{Conj2.7} {\rm (2012-09-08)} Let $f_1(x),\ldots,f_k(x)$ be irreducible polynomials with integer coefficients and positive
leading coefficients such that there is no prime dividing  $\prod_{j=1}^k f_j(q)$ for all $q\in\Z$.
Let $q_1,q_2,\ldots$ be the list (in ascending order) of those  $q\in\Z^+$ such that $f_1(q),\ldots,f_k(q)$ are all primes.
Then, for all sufficiently large positive integers $n$, we have
   $$\f2{n-1}Q(n)<q_{n+1} < q_n^{1+1/n}.$$
Also, for some $N\in\Z^+$  the sequence
  $(\root{n+1}\of{Q(n+1)}/\root n\of{Q(n)})_{n\gs N}$ is strictly increasing with limit $1$.
\end{conjecture}

\begin{remark}\label{Rem2.7} Obviously $2Q(n)<(n-1)q_{n+1}$ if and only if $Q(n+1)<(n+1)q_{n+1}/2$.
\end{remark}

 For convenience, under the condition of Conj. 2.7, below we set
$$E(f_1(x),\ldots,f_k(x))=\{n\in\Z^+:\ \root{n}\of{q_n}>\root{n+1}\of{q_{n+1}}\ \ \t{fails}\}$$
and let $N_0(f_1(x),\ldots,f_k(x))$ stand for the least positive integer $n_0$ such that $2Q(n)<(n-1)q_{n+1}$ for all $n\gs n_0$, and let
$N(f_1(x),\ldots,f_k(x))$ denote the smallest positive integer $N$
such that $(\root{n+1}\of{Q(n+1)}/\root n\of{Q(n)})_{n\gs N}$ is strictly increasing with limit 1.
\medskip

If $p$ and $p+2$ are both primes, then $\{p,p+2\}$ is said to be a pair of twin primes.
The famous twin prime conjecture states that there are infinitely many
twin primes.

\begin{conjecture}\label{Conj2.8} {\rm (2012-08-18)} We have
$$E(x,x+2)=\em,\ \ N_0(x,x+2)=4,\ \ \t{and}\ N(x,x+2)=9.$$
\end{conjecture}

\begin{remark}\label{Rem2.8}  Let  $q_1,q_2,\ldots$ be the list of those primes $p$ with $p+2$ also prime.
We have verified that $\root n\of{q_n}>\root{n+1}\of{q_{n+1}}$ for all $n=1,\ldots,500000$,
$q_{n+1}>2Q(n)/(n-1)$ for all $n=4,\ldots,2000000$, and
$\root{n+1}\of {Q(n+1)}/\root n\of{Q(n)}<\root{n+2}\of {Q(n+2)}/\root {n+1}\of{Q(n+1)}$ for all $n=9,\ldots,500000$.
See also Conjecture 2.10 of the author \cite{S12c}.
\end{remark}

\begin{conjecture}\label{Conj2.9} {\rm (2012-08-20)} We have
\begin{gather*} E(x,x+2,x+6)=E(x,x+4,x+6)=\em,
\\N_0(x,x+2,x+6)=3,\ N_0(x,x+4,x+6)=6,
\\N(x,x+2,x+6)=N(x,x+4,x+6)=13.\end{gather*}
\end{conjecture}

\begin{remark}\label{Rem2.9} Recall that a prime triplet is a set of three primes of the form $\{p,p+2,p+6\}$ or $\{p,p+4,p+6\}$.
It is conjectured that there are infinitely many prime triplets.
\end{remark}

A prime $p$ is called a Sophie Germain prime if $2p+1$ is also a prime.
It is conjectured that there are infinitely many Sophie Germain primes, but this has not
been proved yet.

\begin{conjecture}\label{Conj2.10} {\rm (2012-08-18)}
We have
$$E(x,2x+1)=\{3,4\},\ N_0(x,2x+1)=3,\ \t{and}\ N(x,2x+1)=13.$$
Also,
$$E(x,2x-1)=\{2,3,6\},\ N_0(x,2x-1)=3,\ \t{and}\ N(x,2x-1)=9.$$
\end{conjecture}

\begin{remark}\label{Rem2.10}  When $q_1,q_2,\ldots$ gives the list of Sophie Germain primes in ascending order,
we have verified that $\root n\of{q_n}>\root{n+1}\of{q_{n+1}}$ for all $n=5,\ldots,200000$, and
$\root{n+1}\of {Q(n+1)}/\root n\of{Q(n)}<\root{n+2}\of {Q(n+2)}/\root {n+1}\of{Q(n+1)}$ for every $n=13,\ldots,200000$.
\end{remark}

One may wonder whether $E(x,x+d)$ or $E(x,2x+d)$ with small $d\in\Z^+$ may contain relatively large elements.
We have checked this for $d\ls 100$. Here are few extremal examples suggested by our computation:
\begin{gather*} E(x,x+60)=\{187,3976,58956\},\ E(x,x+66)=\{58616\},
\\E(x,2x+11)=\{1,39593\},\ E(x,2x+81)=\{104260\}.\end{gather*}

\begin{conjecture}\label{Conj2.11} {\rm (2012-09-07)}
We have
\begin{gather*}E(x,x^2+x+1)=\{3,4,12,14\},
\\N_0(x,x^2+x+1)=3,\ N(x,x^2+x+1)=17.
\end{gather*}
Also,
$$E(x^4+1)=\{1,2,4\},\ N_0(x^4+1)=4,\ \t{and}\ N(x^4+1)=10.$$
\end{conjecture}

\begin{remark}\label{Rem2.11} Note that those primes $p$ with $p^2+p+1$ prime are sparser than twin primes and Sophie Germain primes.
\end{remark}

\subsection{Conjectures on other number-theoretic sequences}
$\ $

A positive integer $n$ is called {\it squarefree} if $p^2\nmid n$ for any prime $p$.
Here is the list of all squarefree positive integers not exceeding 30 in ascending order:
$$1,\ 2,\ 3,\ 5,\ 6,\ 7,\ 10,\ 11,\ 13,\ 14,\ 15,\ 17,\ 19,\ 21,\ 22,\ 23,\ 26,\ 29,\ 30.$$

\begin{conjecture}\label{Conj2.12} {\rm (2012-08-14)} Let $s_1,s_2,s_3,\ldots$ be the list of squarefree positive integer in ascending order.
Then the sequence $(\root n\of{s_n})_{n\gs7}$ is strictly decreasing, and
the sequence $(\root{n+1}\of {S(n+1)}/\root n\of{S(n)})_{n\gs7}$ is strictly increasing, where $S(n)=\sum_{k=1}^ns_k$.
\end{conjecture}

\begin{remark}\label{Rem2.12} We have verified that $\root n\of{s_n}>\root{n+1}\of{s_{n+1}}$ for all $n=7,\ldots,500000$.
Note that $\lim_{n\to\infty}\root n\of{S(n)}=1$ since $S(n)$ does not exceed the sum of the first $n$ primes.
\end{remark}

\begin{conjecture}\label{Conj2.13} {\rm (2012-08-25)} Let $a_n$ be the $n$-th positive integer that can be written as a sum of two squares.
Then the sequence $(\root n\of{a_n})_{n\gs6}$ is strictly decreasing, and
the sequence $(\root{n+1}\of {A(n+1)}/\root n\of{A(n)})_{n\gs6}$ is strictly increasing, where $A(n)=\sum_{k=1}^na_k$.
\end{conjecture}

\begin{remark}\label{Rem2.13} Similar things happen if we replace sums of squares in Conj. 2.13 by integers of the form $x^2+dy^2$ with $x,y\in\Z$,
where $d$ is any positive integer.
\end{remark}

Recall that a  partition of a positive integer $n$ is a way of writing $n$ as a sum of  positive integers with the order of addends ignored.
Also, a {\it strict partition} of $n\in\Z^+$ is a way of writing $n$ as a sum of {\it distinct} positive integers with the order of addends ignored.
For $n=1,2,3,\ldots$ we denote by $p(n)$ and $p_*(n)$ the number of partitions of $n$ and the number of strict partitions of $n$ respectively.
It is known that
$$p(n)\sim\f{e^{\pi\sqrt{2n/3}}}{4\sqrt3 n}\ \ \t{and}\ \ p_*(n)\sim\f{e^{\pi\sqrt{n/3}}}{4(3n^3)^{1/4}}\ \quad\t{as}\ n\to+\infty$$
(cf. \cite{HR} and \cite[p.\,826]{AS}) and hence $$\lim_{n\to\infty}\root{n}\of{p(n)}=\lim_{n\to\infty}\root{n}\of{p_*(n)}=1.$$

\begin{conjecture}\label{Conj2.14} {\rm (2012-08-02)} Both $(\root n\of{p(n)})_{n\gs6}$ and $(\root n\of{p_*(n)})_{n\gs9}$
are strictly decreasing. Furthermore, the sequences $(\root{n+1}\of{p(n+1)}/\root n\of{p(n)})_{n\gs26}$ and
$(\root{n+1}\of{p_*(n+1)}/\root n\of{p_*(n)})_{n\gs45}$ are strictly increasing.
\end{conjecture}

\begin{remark}\label{Rem2.14} The author has verified the conjecture for $n$ up to $10^5$. \cite{S12c} contains a stronger version of this conjecture.
\end{remark}

The Bernoulli numbers $B_0,B_1,B_2,\ldots$ are rational numbers given by
$$B_0=1,\ \ \t{and}\ \ \sum^n_{k=0}\bi {n+1}k B_k=0\ \ \ \t{for}\ n\in\Z^+.$$
It is well known that $B_{2n+1}=0$ for all $n\in\Z^+$ and
$$\f x{e^x-1}=\sum_{n=0}^\infty B_n\f{x^{n}}{n!}\ \ \l(|x|<2\pi\r).$$
(See, e.g., \cite[pp.\,228-232]{IR}.)
The Euler numbers $E_0,E_1,E_2,\ldots$ are integers defined by
$$E_0=1,\ \ \t{and}\ \ \sum^n_{k=0\atop 2\mid k} \bi nk E_{n-k}=0\ \ \ \t{for}\ n\in\Z^+.$$
It is well known that $E_{2n+1}=0$ for all $n=0,1,2,\ldots$ and
$$\sec x=\sum_{n=0}^\infty(-1)^n E_{2n}\f{x^{2n}}{(2n)!}\ \ \l(|x|<\f{\pi}2\r).$$

\begin{conjecture}\label{Conj2.15} {\rm (2012-08-02)} $(\root n\of{(-1)^{n-1}B_{2n}})_{n\gs1}$ and $\root n\of{(-1)^nE_{2n}})_{n\gs1}$
are strictly increasing, where $B_0,B_1,\ldots$ are Bernoulli numbers
and $E_0,E_1,\ldots$ are Euler numbers. Moreover, the sequences
$$\l(\root{n+1}\of{(-1)^{n}B_{2n+2}}/\root n\of{(-1)^{n-1}B_{2n}}\r)_{n\gs2}$$
and $$\l(\root{n+1}\of{(-1)^{n+1}E_{2n+2}}/\root n\of{(-1)^nE_{2n}}\r)_{n\gs1}$$ are strictly decreasing.
\end{conjecture}

\begin{remark}\label{Rem2.15} It is known that both $(-1)^{n-1}B_{2n}$ and $(-1)^nE_{2n}$ are positive for all $n=1,2,3,\ldots$.
\end{remark}

For $m,n\in\Z^+$ the $n$-th harmonic number $H_n^{(m)}$ of order $m$ is defined as $\sum_{k=1}^n1/k^m$.

\begin{conjecture}\label{Conj2.16} {\rm (2012-08-12)} For any positive integer $m$, the sequence
$$\l(\root{n+1}\of{H_{n+1}^{(m)}}\big/\root{n}\of{H_n^{(m)}}\r)_{n\gs3}$$ is strictly increasing.
\end{conjecture}

\begin{remark}\label{Rem2.16} It is easy to show that $\big(\root n\of{H_n^{(m)}}\big)_{n\gs2}$ is strictly decreasing for any $m\in\Z^+$.
Some fundamental congruences on harmonic numbers can be found in \cite{S12a}.
\end{remark}

\begin{conjecture}\label{Conj2.17} {\rm (2012-09-01)} Let $q>1$ be a prime power and let $\Bbb F_q$ be the finite field of order $q$.
Let $M_n(q)$ denote the number of monic irreducible polynomials of degree at most $n$ over $\Bbb F_q$.

{\rm (i)} We have $ M_q(n+1)/M_q(n) < M_q(n+2)/M_q(n+1)$
unless $q<5$ and $n\in\{2,4,6,8,10,12\}$.

{\rm (ii)} If $n>2$, then
$\root n\of{M_q(n)}< \root{n+1}\of{M_q(n+1)}$
unless $q<7$ and $n\in\{3,5\}$.

{\rm (iii)} When $n>3$, we have
   $$\root{n+1}\of{M_q(n+1)}\big/\root n\of{M_q(n)}
  > \root{n+2}\of{M_q(n+2)}\big/\root{n+1}\of{M_q(n+1)}$$
unless $(q<8\ \&\ n\in\{5,7,9,11,13\})$ or $(9<q<14\ \&\ n=4)$.
\end{conjecture}

\begin{remark}\label{Rem2.17} It is known that the number of monic irreducible polynomials of degree $n$ over the finite field $\Bbb F_q$
equals $\f1n\sum_{d\mid n}\mu(d)q^{n/d}$, where $\mu$ is the M\"obius function (cf. \cite[p.\,84]{IR}).
\end{remark}

\section{Conjectures on combinatorial sequences}
\setcounter{lemma}{0}
\setcounter{theorem}{0}
\setcounter{corollary}{0}
\setcounter{remark}{0}
\setcounter{equation}{0}
\setcounter{conjecture}{0}

The Fibonacci sequence $(F_n)_{n\gs0}$ is given by
$$F_0=0,\ F_1=1,\ \t{and}\ F_{n+1}=F_n+F_{n-1}\ (n=1,2,3,\ldots).$$
the reader may consult \cite[p.\,46]{St} for combinatorial interpretations of Fibonacci
numbers.

\begin{conjecture}\label{Conj3.1} {\rm (2012-08-11)} The sequence
$(\root{n}\of{F_n})_{n\gs2}$ is strictly increasing, and moreover the
sequence $(\root{n+1}\of {F_{n+1}}/\root n\of{F_n})_{n\gs 4}$ is
strictly decreasing. Also, for any integers $A>1$ and $B\not=0$ with $A^2>4B$ and $(A>2$ or $B\gs-9)$,
the
sequence $(\root{n+1}\of {u_{n+1}}/\root n\of{u_n})_{n\gs 4}$ is
strictly decreasing with limit $1$, where
$$u_0=0,\ u_1=1,\ \t{and}\ u_{n+1}=Au_n-Bu_{n-1}\ (n=1,2,3,\ldots).$$
\end{conjecture}

\begin{remark}\label{Rem3.1} By \cite[Lemma 4]{S92}, if $A>1$ and $B\not=0$ are integers with $A^2>4B$
then the sequence $(u_n)_{n\gs0}$ defined in Conjecture 3.1
is strictly increasing.
\end{remark}

For $n=1,2,3,\ldots$ the $n$-th Bell number $B_n$ denotes the number
of partitions of $\{1,\ldots,n\}$ into disjoint nonempty subsets.
It is known that $B_{n+1}=\sum_{k=0}^n\bi nk B_k$ (with $B_0=1$) and
$B_n=e^{-1}\sum_{k=0}^\infty k^n/k!$ for all $n=0,1,2,\ldots$ (cf. \cite[A000110]{Sl}).

\begin{conjecture}\label{Conj3.2} {\rm (2012-08-11)} The sequence
$(\root{n}\of{B_n})_{n\gs1}$ is strictly increasing, and moreover
the sequence $(\root{n+1}\of {B_{n+1}}/\root n\of{B_n})_{n\gs 1}$ is
strictly decreasing with limit $1$, where $B_n$ is the $n$-th Bell number.
\end{conjecture}

\begin{remark}\label{Rem3.2} In 1994 K. Engel \cite{E} proved the log-convexity of $(B_n)_{n\gs1}$.
\cite{SZ} contains a curious congruence property of the Bell numbers.
\end{remark}

For $n\in\Z^+$ the $n$-th derangement number $D_n$ denotes the number
of permutations $\sigma$ of $\{1,\ldots,n\}$ with $\sigma(i)=i$ for
no $i=1,\ldots,n$. It has the following explicit expression (cf. \cite[p.\,67]{St}):
$$D_n=\sum_{k=0}^n(-1)^k\f{n!}{k!}.$$

\begin{conjecture}\label{Conj3.3} {\rm (2012-08-11)} The sequence
$(\root{n}\of{D_n})_{n\gs2}$ is strictly increasing, and the
sequence $(\root{n+1}\of {D_{n+1}}/\root n\of{D_n})_{n\gs 3}$ is
strictly decreasing.
\end{conjecture}

\begin{remark}\label{Rem3.3} As $D_n=nD_{n-1}+(-1)^n$ for $n\in\Z^+$, it is easy to see that $(D_{n+1}/D_n)_{n\gs1}$ is strictly increasing.
\end{remark}

During his study of irreducible root systems of a special type related to Weyl groups,  T. A. Springer \cite{Sp}
introduced the Springer numbers $S_0,S_1,\ldots$ defined by
$$\f1{\cos x-\sin x}=\sum_{n=0}^\infty S_n\f{x^n}{n!}.$$
The reader may consult \cite[A001586]{Sl} for various combinatorial interpretations of Springer numbers.

\begin{conjecture}\label{Conj3.4} {\rm (2012-08-05)} The sequence
$(S_{n+1}/S_n)_{n\gs0}$ is strictly increasing, and the sequence
$(\root{n+1}\of {S_{n+1}}/\root n\of{S_n})_{n\gs 1}$ is strictly
decreasing with limit $1$, where $S_n$ is the $n$-th Springer number.
\end{conjecture}

\begin{remark}\label{Rem3.4} It is known (cf. \cite[A001586]{Sl}) that $S_n$ coincides with the numerator of $|E_n(1/4)|$,
where $E_n(x)$ is the Euler polynomial of degree $n$.
\end{remark}

\begin{conjecture}\label{Conj3.5} {\rm (2012-08-18)} For the tangent numbers $T(1),T(2),\ldots$ given by
$$\tan x=\sum_{n=1}^\infty T(n)\f{x^{2n-1}}{(2n-1)!},$$
the sequences $(T(n+1)/T(n))_{n\gs1}$ and  $(\root n\of{T(n)})_{n\gs1}$ are strictly increasing,
 and the sequence $(\root{n+1}\of{T(n+1)}/\root n\of{T(n)})_{n\gs 2}$
is strictly decreasing.
\end{conjecture}

\begin{remark}\label{Rem3.5} The tangent numbers are all integral, see \cite[A000182]{Sl} for the sequence $(T(n))_{n\gs1}$.
It is known that $T(n)=(-1)^{n-1}2^{2n}(2^{2n}-1)B_{2n}/(2n)$ for all $n\in\Z^+$, where $B_{2n}$ is the $2n$-th Bernoulli number.
\end{remark}

The $n$-th central trinomial coefficient $T_n$ is the coefficient of
$x^n$ in the expansion of $(x^2+x+1)^n$. Here is an explicit
expression:
$$T_n=\sum_{k=0}^n\bi nk\bi{n-k}k=\sum_{k=0}^{\lfloor n/2\rfloor}\bi
n{2k}\bi{2k}k.$$ In combinatorics, $T_n$ is the number of lattice paths from the
point $(0, 0)$ to $(n, 0)$ with only allowed steps $(1,0)$, $(1, 1)$
and $(1, -1)$ (cf. \cite[A002426]{Sl}). It is known that $(n+1)T_{n+1}=(2n+1)T_n+3nT_{n-1}$ for all $n\in\Z^+$.

\begin{conjecture}\label{Conj3.6} {\rm (2012-08-11)} The sequence
$(\root{n}\of{T_n})_{n\gs1}$ is strictly increasing, and
the sequence $(\root{n+1}\of {T_{n+1}}/\root n\of{T_n})_{n\gs 1}$ is
strictly decreasing.
\end{conjecture}

\begin{remark}\label{Rem3.6} Via the Laplace-Heine formula (cf. \cite[p.\,194]{Sz}) for
Legendre polynomials, $T_n\sim3^{n+1/2}/(2\sqrt{n\pi})$ as
$n\to+\infty$. In 2011, the author \cite{S11c} found many series for
$1/\pi$ involving generalized central trinomial coefficients.
\end{remark}

The $n$-th Motzkin number
$$M_n=\sum_{k=0}^{\lfloor n/2\rfloor}\bi n{2k}\bi{2k}k\f1{k+1}$$
is the number of lattice paths from $(0,0)$ to $(n,0)$ which never
dip below the line $y=0$ and are made up only of the allowed steps
$(1,0)$, $(1,1)$ and $(1,-1)$ (cf. \cite[A001006]{Sl}). It is known that $(n+3)M_{n+1}=(2n+3)M_n+3nM_{n-1}$ for all $n\in\Z^+$.

\begin{conjecture}\label{Conj3.7} {\rm (2012-08-11)} The sequence
$(\root{n}\of{M_n})_{n\gs1}$ is strictly increasing, and moreover
the sequence $(\root{n+1}\of {M_{n+1}}/\root n\of{M_n})_{n\gs 1}$ is
strictly decreasing.
\end{conjecture}

\begin{remark}\label{Rem3.7} The log-convexity of the sequence $(M_n)_{n\gs1}$ was first established by M. Aigner \cite{A} in 1998.
\end{remark}

For $r=2,3,4,\ldots$ define
$$f_n^{(r)}:=\sum_{k=0}^n\bi nk^r\ \ (n=0,1,2,\ldots).$$
Note that $f_n^{(2)}=\bi{2n}n$, and those $f_n=f_n^{(3)}$ are called
Franel numbers (cf. \cite[A000172]{Sl}).

\begin{conjecture}\label{Conj3.8} {\rm (2012-08-11)} For each $r=2,3,4,\ldots$
there is a positive integer $N(r)$ such that the sequence
$\big(\root{n+1}\of {f_{n+1}^{(r)}}/\root n\of{f_n^{(r)}}\big)_{n\gs
N(r)}$ is strictly decreasing with limit $1$. Moreover, we may take
\begin{gather*} N(2)=\cdots=N(6)=1,\ \  N(7)=N(8)=N(9)=3,\ N(10)=N(11)=5,
\\ N(12)=N(13)=7,\ N(14)=N(15)=N(16)=9,\ N(17)=N(18)=11.
\end{gather*}
\end{conjecture}

\begin{remark}\label{Rem3.8} It is known that $(f_n^{(r)})_{n\gs1}$ is log-convex for $r=2,3,4$ (cf. \cite{D05}).
\cite{S11b} contains some fundamental congruences
for Franel numbers.
\end{remark}

\begin{conjecture}\label{Conj3.9} {\rm (2012-08-15)} Set $g_n=\sum_{k=0}^n\bi nk^2\bi{2k}k$ for $n=0,1,2,\ldots$.
Then  $(\root n\of{g_n})_{n\gs1}$
is strictly increasing and the sequence $(\root{n+1}\of{g_{n+1}}/\root{n}\of{g_n})_{n\gs1}$ is strictly decreasing.
\end{conjecture}

\begin{remark}\label{Rem3.9} It is known that $g_n=\sum_{k=0}^n\bi nkf_k$, where $f_k=\sum_{j=0}^k\bi kj^3$ is the $k$-th Franel number.
Both $(f_n)_{n\gs0}$ and $(g_n)_{n\gs0}$ are related to the theory of modular
forms, see D. Zagier \cite{Z}.
\end{remark}

For $r=1,2,3,\ldots$ define
$$A_n^{(r)}=\sum_{k=0}^n\bi nk^r\bi{n+k}k^r\ \ (n=0,1,2,\ldots).$$
Those $A_n^{(1)}$ and $A_n=A_n^{(2)}$ are called central Delannoy
numbers and Ap\'ery numbers respectively. The Ap\'ery numbers play a key role in Ap\'ery's proof of the irrationality of
$\zeta(3)=\sum_{n=1}^\infty1/n^3$ (cf. \cite{Ap,P}).

\begin{conjecture}\label{Conj3.10} {\rm (2012-08-11)} For each $r=1,2,3,\ldots$
there is a positive integer $M(r)$ such that the sequence
$\big(\root{n+1}\of {A_{n+1}^{(r)}}/\root n\of{A_n^{(r)}}\big)_{n\gs
M(r)}$ is strictly decreasing with limit $1$. Moreover, we may take
$$M(1)=\cdots=M(16)=1,\ \  M(17)=M(18)=M(19)=9,\ M(20)=12.$$
\end{conjecture}

\begin{remark}\label{Rem3.10} The log-convexity of $(A_n)_{n\gs0}$ was proved by T. Do\v sli\'c \cite{D05}.
The reader may consult \cite{S12b} for some congruences involving Ap\'ery numbers and Ap\'ery polynomials.
\end{remark}

The $n$-th Schr\"oder number
$$S_n=\sum_{k=0}^n\bi
nk\bi{n+k}k\f1{k+1}=\sum_{k=0}^n\bi{n+k}{2k}\bi{2k}k\f1{k+1}$$ is
the number of lattice paths from the point $(0,0)$ to $(n,n)$ with steps $(1,0),(0,1)$ and $(1,1)$ that never rise above the line $y=x$
(cf. \cite[A006318]{Sl} and \cite[p.\,185]{St}).

\begin{conjecture}\label{Conj3.11} {\rm (2012-08-11)} The sequence
$(\root{n}\of{S_n})_{n\gs1}$ is strictly increasing, and moreover
the sequence $(\root{n+1}\of {S_{n+1}}/\root n\of{S_n})_{n\gs 1}$ is
strictly decreasing, where $S_n$ stands for the $n$-th Schr\"oder number.
\end{conjecture}

\begin{remark}\label{Rem3.11} The reader may consult \cite{S11a} for some congruences involving central Delannoy numbers and Schr\"oder numbers.
\end{remark}

\begin{conjecture}\label{Conj3.12} {\rm (2012-08-13)} For the Domb numbers
$$D(n)=\sum_{k=0}^n\bi nk^2\bi{2k}k\bi{2(n-k)}{n-k}\ (n=0,1,2,\ldots),$$
the sequences $(D(n+1)/D(n))_{n\gs0}$ and $(\root n\of{D(n)})_{n\gs1}$ are strictly increasing. Moreover,
the sequence $(\root{n+1}\of{D(n+1)}/\root n\of {D(n)})_{n\gs1}$ is strictly decreasing.
\end{conjecture}

\begin{remark}\label{Rem3.12} For combinatorial interpretations of the Domb number $D(n)$, the reader may consult \cite[A002895]{Sl}.
\cite{CCL} contains some series for $1/\pi$ involving Domb numbers.
\end{remark}

The Catalan-Larcombe-French numbers $P_0,P_1,P_2,\ldots$ (cf. \cite{JV}) are given by
$$P_n=\sum_{k=0}^n\f{\bi{2k}k^2\bi{2(n-k)}{n-k}^2}{\bi nk}=2^n\sum_{k=0}^{\lfloor n/2\rfloor}\bi n{2k}\bi{2k}k^24^{n-2k},$$
they arose from the theory of elliptic integrals (see \cite{LF}).
It is known that $(n+1)P_{n+1}=(24n(n+1)+8)P_n-128n^2P_{n-1}$ for all $n\in\Z^+$.
The sequence $(P_n)_{n\gs0}$ is also related to the theory of modular forms, see D. Zagier \cite{Z}.

\begin{conjecture}\label{Conj3.13} {\rm (2012-08-14)} The sequences $(P_{n+1}/P_n)_{n\gs0}$
and $(\root n\of{P_n})_{n\gs1}$ are strictly increasing. Moreover,
the sequence $(\root{n+1}\of{P_{n+1}}/\root n\of {P_n})_{n\gs1}$ is strictly decreasing.
\end{conjecture}

\begin{remark}\label{Rem3.13}
We also have the following conjecture related to Euler numbers:
$$\sum_{k=0}^{p-1}\f {P_k}{8^k}\eq1+2\l(\f{-1}p\r)p^2E_{p-3} \ (\mo\ p^3)$$
and
$$\sum_{k=0}^{p-1}\f {P_k}{16^k}\eq\l(\f{-1}p\r)-p^2E_{p-3} \ (\mo\
p^3)$$ for any odd prime $p$, where $(\f{\cdot}p)$ is the Legendre symbol.
\end{remark}

\Ack. The initial work was done during the author's visit to the
University of Illinois at Urbana-Champaign, so the author wishes to
thank Prof. Bruce Berndt for his kind invitation and hospitality.
The author is also grateful to the referee for helpful comments.

\end{document}